\documentclass[a4paper,10pt]{article}
\usepackage{amsmath,amscd,amsthm,amsfonts,amssymb,epsf,latexsym}
\usepackage[all]{xy}

\makeatletter

\long\def\@makefnt#1{\parindent 1em\noindent
             \hb@xt@1.8em{\hss\@textsuperscript{}}#1}
\long\def\@ftntext#1{\insert\footins{%
     \reset@font\footnotesize
     \interlinepenalty\interfootnotelinepenalty
     \splittopskip\footnotesep
     \splitmaxdepth \dp\strutbox \floatingpenalty \@MM
     \hsize\columnwidth \@parboxrestore
     \color@begingroup
       \@makefnt{%
         \rule\z@\footnotesep\ignorespaces#1\@finalstrut\strutbox}%
     \color@endgroup}}%
\def\subjclass#1{%
   \@ftntext{2010 {\itshape Mathematics Subject Classification.}\enspace
#1.}}
\def\keywords#1{%
   \@ftntext{{\itshape Keywords.}\enspace #1.}}
\makeatother

\def\moins{\raise 1pt\hbox{{$\scriptstyle -$}}}
\def\plus{\raise 1pt\hbox{{$\scriptstyle +$}} }

\newtheorem{theorem}{Theorem}
\newtheorem{proposition}[theorem]{Proposition}
\newtheorem{lemma}[theorem]{Lemma}
\newtheorem{corollary}[theorem]{Corollary}
\newtheorem{remark}[theorem]{Remark}

\newtheorem{example}[theorem]{Example}

\def\proof{\noindent{\bf Proof.\ }}

\def\qed{~\hbox{$\Box$}}

\def\Aut{\mathop{\rm Aut}}
\def\cJ{\mathop{\rm J}}
\def\codim{\mathop{\rm codim}}

\def\Ker{\mathop{\rm Ker}}
\def\dim{\mathop{\rm dim}}
\def\rank{\mathop{\rm rank}}
\def\card{\mathop{\rm card}}

\def\Sympl{\mathop{\rm Sympl}}
\def\End{\mathop{\rm End}}

\let\L\l
\def\l{\mathop{\lambda}}

\def\cT{{\mathcal T}}
\def\cJ{{\mathcal J}}

\def\cS{{\mathcal S}}
\def\cO{{\mathcal O}}
\def\cF{{\mathcal F}}

\def\C{{\bf C}}
\def\Z{{\bf Z}}
\def\N{{\bf N}}
\def\R{{\bf R}}
\def\P{{\bf P}}
\def\Q{{\bf Q}}
\def\G{{\bf G}}

\def\\Q{\widetilde{Q}}

\def\qed{~\hbox{$\Box$}}

\begin{document}

\title{\bf Positivity of Thom polynomials and Schubert calculus}

\author{
Piotr Pragacz\thanks{Research supported by a MNiSzW grant N N201 608040.}\\
\small Institute of Mathematics, Polish Academy of Sciences\\
\small \'Sniadeckich 8, 00-656 Warszawa, Poland\\
\small P.Pragacz@impan.pl}



\date{}

\maketitle

%\tableofcontents

\begin{abstract}
We describe the positivity of Thom polynomials of singularities of maps, Lagrangian Thom polynomials 
and Legendrian Thom polynomials. We show that these positivities come from Schubert calculus. 

\end{abstract}

\section{Introduction}

In the present paper, we discuss the issue of {\it positivity}. The positivity plays an important
role in mathematics. For example, positivity in algebraic geometry is a subject of a vast monograph \cite{Lz} 
of Lazarsfeld. 

There are two questions related to positivity: 1. Are the numbers in question (mostly the coefficients of some polynomials) nonnegative?
2. If yes, what is a positive description of these numbers?

Positivity in Schubert calculus is an active area of the contemporary research, related mainly to the second question,
see, e.g., the survey article \cite{CoVa}. Answers to question 1. in many important situations
are known classically (many of them follow from the Bertini-Kleiman theorem). Since the author hopes that the present
article will also be read by beginners, we discuss this issue briefly at the end of Section \ref{Schub}.

Our main goal here is to describe some positivities in the global geometry of singularities. 
These positivities come from Schubert calculus. The presented results are related mainly to the first question.
The present knowledge about the second question in this area is rather restricted.
We survey several recent positivity results about Thom polynomials.
Some of them are obtained by the Bertini-Kleiman theorem and its variants; other are deduced using the Fulton-Lazarsfeld 
theorems on positive polynomials for ample vector bundles, or vector bundles generated by their global sections. 
We also discuss positivity of the restrictions of Schubert classes, and some other related positivities. 

Here is a description of the content of the paper.

After the preliminaries, in Section \ref{Schub}, we recall basic definitions and facts about Schubert classes in the cohomology rings
of $G/P$. We put an emphasis on Poincar\'e duality.
Then -- to start the discussion on positivity -- we recall two known positivities in the cohomology rings of flag manifolds
and explain the positivity of restrictions of Schubert classes.

Thom polynomials came from algebraic topology and singularities. The classical ones are associated with singularities of maps.
Nowadays, we also study Lagrangian Thom polynomials and Legendrian Thom polynomials. In Sections \ref{tsm}, \ref{tlag} and \ref{tleg}, we
give brief introductions to these three series of Thom polynomials. The computations of Thom polynomials are, in general, 
quite hard. There are basically two ways to compute the Thom polynomials of a singularity class $\Sigma$: 1. using desingularization
of $\Sigma$, and push-forward formulas (one should mention here many names, making this article too long); 2. the interpolation method
of Feh\'er and Rim\'anyi: by restricting to singularities of smaller codimension than $\codim \Sigma$, and using symmetries of singularities (see \cite{Rim}). 
It was the basis of monomials in Chern classes, which served at first to compute the Thom
polynomials. About a decade ago, the basis of Schur functions started also to be used systematically for computations
of Thom polynomials (see Section \ref{tsm} for more details). 

In 2006 Weber and the author proved the positivity of the Thom polynomials of stable singularity classes of maps in the basis of Schur functions
\cite{PrWe1}. The method relies on classifying spaces of singularities and on some global aspects of Schubert calculus.
The Fulton-Lazarsfeld theory \cite{FuLa1} of polynomials numerically positive for {\it ample} vector bundles is used. For details, see Sections \ref{amppos} and \ref{tsm}. Thus methods of algebraic geometry appear to be useful to study Thom polynomials.

Section \ref{tic} presents a generalization, by the same two co-authors, of this positivity to Thom polynomials of invariant cones 
and, in particular, to the Thom polynomials of possibly nonstable singularity classes of maps \cite{PrWe2}.

In section \ref{tlag}, we describe the positivity of Lagrangian Thom polynomials in the basis of $\\Q$-polynomials.
This is a result of Mikosz, Weber and the author \cite{MPW1}.

The positivity of Legendrian Thom polynomials is a subject of Section \ref{tleg}, where we report on results of the same three 
co-authors \cite{MPW2}. The argument is based on some variant of the Bertini-Kleiman theorem and the Schubert calculus for Lagrangian 
Grassmann bundle associated with a twisted skew-symmetric form. One constructs a basis in the cohomology ring
of that Lagrangian Grassmann bundle such that any Legendrian Thom polynomial has, in this basis, an expansion with nonnegative 
coefficients. This leads to the construction of a one-parameter family of such bases in the ring of Legendrian characteristic classes.

\smallskip

This is a written account of the talk delivered by the author at the conference on Schubert calculus (July 23-27, 2012)
at Osaka, where apart from interesting lectures he enjoyed beautiful Japanese gardens (cf. \cite{Shi}).
He wishes to thank the organizers of the conference for their devoted work.

\smallskip

The main body of the paper is based on a cooperation with Ma\L gorzata Mikosz and Andrzej Weber, with some assistance of Maxim
Kazarian and Alain Lascoux. The author is grateful to them for invaluable conversations. He also thanks Wojciech Domitrz, Letterio Gatto, 
Megumi Harada, Jaros\L aw K\c edra, Toshiaki Maeno, Piotr Mormul and Krzysztof Pawa\L owski for useful comments. 
Finally, the author thanks the referee for pointing out several defects in the previous version of this article and suggesting
some improvements.

\section{Preliminaries}

General information about varieties, homology groups $H_*(-)$, cohomology groups $H^*(-)$ and Chow groups $A_*(-)$,
$A^*(-)$ 
in the scope needed for this paper is contained in \cite[App. A]{FuPr}. The multiplication in cohomology and that in Chow rings
of nonsingular varieties will be denoted by ``$\cdot$''. For more detailed information concerning these matters, we refer
the reader to \cite{Fu1}. We follow the notation for algebraic geometry from this book. We use the following variants of fundamental classes:

\begin{itemize}

\item Let $X$ be a variety over a field $k$. Given a (closed) subscheme $Z$ of $X$ of pure dimension $d$, by $[Z]=\sum m_i[Z_i]$
we denote its fundamental class in the Chow group $A_d(X)$, where $Z_i$ are irreducible components of $Z$ and
$m_i=l(\cO_{Z,Z_i})$ are their geometric multiplicities.

\item If $k=\C$, a (closed) subscheme $Z$ of a compact variety $X$ of pure dimension $d$ determines in the same way a fundamental class
in $H_{2d}(X,\Z)$ denoted $[Z]$. If $X$ is nonsingular, by Poincar\'e duality, we have the class $[Z]$ in
$H^{2e}(X,\Z)=H_{2d}(X,\Z)$, where $e$ is the codimension of $Z$ in $X$.

\smallskip

\item If $Z$ is a (closed) subscheme of a possibly {\it noncompact} complex manifold\footnote{A manifold here is always {\it nonsingular}.} $X$
of pure codimension $e$, then we have a class 
$[Z] \in H^{2e}(X,\Z)$. Indeed, $Z$ has a fundamental class $[Z]$ in the Borel-Moore homology group $H_{2d}(X)$,
$d=\dim(Z)$ (see \cite{BoHa}), and that group is naturally isomorphic to $H^{2e}(X,\Z)$ (see \cite[Thm 7.9]{BoMo}). 

\end{itemize}
\smallskip

A cycle $\sum n_i [V_i]$ on a scheme $X$ is {\it nonnegative} if each $n_i$ is nonnegative.

\medskip

Let $E$ and $F$ be vector bundles on a nonsingular variety $X$.
We define two families of symmetric functions: $s_\lambda(E-F)\in A^{|\lambda|}(X)$ and $\\Q_\mu(E)\in A^{|\mu|}(X)$. 
We follow the notation for partitions from \cite{FuPr}.
Let $\{e\}$ and $\{f \}$ be the sets of Chern roots of $E$ and $F$. We set
\begin{equation}
\sum s_i(E-F) z^i:=\prod_{f} (1-fz)/\prod_{e}(1-ez)\,,
\end{equation}
where $z$ is an independent variable. 
We see that $s_i(E-F)$ interpolates between $s_i(E)$ -- the $i$-th Segre class of $E$ times $(-1)^i$ (cf. \cite{Fu1})
and $s_i(-F)$ -- the $i$-th Chern class of $F$ times $(-1)^i$ ({\it loc.cit.}).
Given a partition 
$$
\lambda=(\lambda_1 \ge \lambda_2\ge \cdots \ge \lambda_l\ge 0)\,,
$$ 
we define
\begin{equation}\label{schur}
s_\lambda(E - F):= \bigl| s_{\lambda_{i}-i+j}(E-F) \bigr|_{1\leq i,j \le l} \, .
\end{equation}
If $F=0$, we obtain $s_{\lambda}(E)$, that is, the classical Schur function of $E$. In the following, the reader will find
a formula, how, knowing the Chern classes of $E$, to get $s_\lambda(E)$. For more detailed information
on the {\it supersymmetric} Schur functions $s_\lambda(E-F)$, see \cite[Sect. 3, 4 and 5]{FuPr}. 

\medskip
We now define the second family of functions.
We set $\widetilde Q_i(E)=c_i(E)$. Given two nonnegative integers $i\ge j$, we define 
$$
\widetilde Q_{i,j}(E):=\widetilde Q_i(E)\widetilde Q_j(E) + 2\sum\limits^j_{p=1}(-1)^p\widetilde Q_{i+p}(E)\widetilde Q_{j-p}(E)\,.
$$
For a partition $\mu$, $\widetilde Q_\mu(E)$ is defined recurrently on $l(\mu)$, 
by putting for odd $l(\mu)$,
$$
\widetilde Q_\mu(E) = \sum\limits_{p=1}^{l(\mu)}(-1)^{p-1} \widetilde Q_{\mu_p}(E) \ \widetilde Q_{\mu \smallsetminus \{\mu_p\}}(E),
$$
and for even $l(\mu)$,
$$
\widetilde Q_\mu(E) = \sum\limits_{p=2}^{l(\mu)} (-1)^p \widetilde Q_{\mu_1,\mu_p}(E) \ \widetilde Q_{\mu \smallsetminus\{\mu_1,\mu_p\}}(E). 
$$
This family of functions is modeled on Schur $Q$-functions, and is useful in Schubert calculus of Lagrangian Grassmannians. The reader can find in \cite[Sect. 3 and 7]{FuPr} more details concerning the polynomials $\widetilde Q_\mu(E)$.

\section{Schubert varieties and Schubert classes}\label{Schub}

In this section, we collect basic information on the cohomology rings of the flag manifolds $G/P$.
We begin by fixing some notation. 

Let $G$ be a semisimple group over an algebraically closed field $k$, and $B\subset G$ a Borel subgroup. Choose a maximal torus $T\subset B$
with Weyl group $W=N_G(T)/T$ of $(G,T)$. This determines a root system $R$, simple roots $\Delta$, positive roots $R^+$ etc. 
The group $W$ is generated by simple reflections $\{s_\alpha: \alpha\in \Delta \}$ with length function $l(w)$ and longest 
word $w_0$: $l(w_0)=\card(R^+)$. The Chevalley-Bruhat decomposition $G=BWB$ provides a ``cell-decomposition'' of the flag manifold
$$
G/B=\coprod_{w\in W} BwB/B\,.
$$
Recall that the flag manifold $G/B$ is nonsingular algebraic and projective of dimension $\card(R^+)$. 
Each subset $Bw_0wB/B$ of $G/B$ is isomorphic, as a $k$-variety, to the affine space \ $k^{l(w_0)-l(w)}$. Its closure 
$\overline{Bw_0wB/B}$ is called a {\it Schubert variety}. This is, in general, a singular algebraic variety of codimension $l(w)$
in $G/B$.  We set in $A^{l(w)}(G/B)$, or in $H^{2l(w)}(G/B,\Z)$
$$
X^w:=[\overline{Bw_0wB/B}]\,,
$$
and call it a {\it Schubert class}. 

The same applies to all partial flag manifolds $G/P$, where $P$ is a parabolic subgroup of $G$. 
Let $\theta$ be a subset of $\Delta$ and let $W_\theta$ be the subgroup of $W$ generated by $\{s_\alpha\}_{\alpha\in\theta}$.
We set $P=P_\theta = BW_\theta B$, and $W_P=W_\theta$.

Consider the set
$$
W^P=W^\theta:= \{w\in W: \ l(ws_\alpha)=l(w)+1\quad \forall\alpha\in\theta\}.
$$
This is the set of minimal length left coset representatives of $W_P$ in $W$.
The projection $G/B\to G/P$ induces an injection
$$
A^{*}(G/P)\hookrightarrow A^{*}(G/B)\,,
$$
which additively identifies $A^{*}(G/P)$ with $\bigoplus_{w\in W^P} {\bf Z} X^w$.
In other words, the $X^w$, $w\in W^P$, form a $\Z$-basis for $A^*(G/P)$ \cite[Thm 5.5]{BGG}.
Multiplicatively, $A^{*}(G/P)_\Q$ is identified with the ring of invariants ${A^{*}(G/B)_\Q}^{W_P}$
\cite[Sect. 5]{BGG}.

\smallskip

If $k=\C$, since $G/P$ admits a cell-decomposition, we have
\begin{equation}\label{ha}
H^{2i+1}(G/P,\Z)=0 \ \ \ \ \hbox{and} \ \ \ \ H^{2i}(G/P,\Z)=A^i(G/P)
\end{equation}
(cf. \cite[Ex. 19.1.11]{Fu1}).

\begin{example}\label{Gr} \rm
Let $G=SL_n$.
We set $P=P_{\theta}$, where $\theta$ is obtained by omitting
the root $\varepsilon_r-\varepsilon_{r+1}$ in the basis
$\varepsilon_1-\varepsilon_2,\ldots,
\varepsilon_{n-1}-\varepsilon_n$ of the root system of type $A_{n-1}$:
$$
\{\varepsilon_i-\varepsilon_j \ | \ i \ne j \}
\subset \oplus_{i=1}^n {\R}\varepsilon_i\,.
$$
We have an identification $SL_n/P=G_r(k^n)$, the Grassmannian parametrizing $r$-dimensional linear subspaces in $k^n$.
It is an algebraic variety of dimension $r(n-r)$. 

The Weyl group $W$ is here the symmetric group $S_n$, and $W_P=S_r\times S_{n-r}$.
The poset $W^P$ is naturally identified with the poset of partitions $\l$ contained in $((n-r)^r)$ (see, e.g., \cite{Hi})
and the corresponding Schubert class $X^{\lambda}$ is represented by the following locus in the Grassmannian.
Consider a flag
$$
V_0\subset V_1 \subset \cdots \subset V_n=k^n
$$
of vector spaces with $\dim(V_i)=i$. Consider the following locus:
$$
\{L\in G_r(k^n) : \dim(L\cap V_{n-r+i-\lambda_i}) \ge i,  \ 1\le i \le r\}\,.
$$
The class of this locus does not depend on the flag, and is equal to the Schubert class $X^\lambda$.

\begin{theorem} [Giambelli formula]\label{Gia} In $A^{|\lambda|}(G_r(k^n))$, we have
$$
X^\lambda =s_\lambda(R^*)\,,
$$
where $R$ is the tautological subbundle on the Grassmannian.
\end{theorem}

\end{example}

\begin{example}\label{LG} \rm
Let $V$ be a symplectic vector space over $k$ of dimension $2n$, and let $G=Sp(V)$ be the symplectic group.
We set $P=P_{\theta}$, where $\theta$ is obtained by omitting
the root $2\varepsilon_n$ in the basis
$\varepsilon_1-\varepsilon_2,\ldots,
\varepsilon_{n-1}-\varepsilon_n, 2\varepsilon_n$ of the root system of type $C_{n}$:
$$
\{\pm \varepsilon_i\pm \varepsilon_j : 1\le i\le j\le n\} \cap \{\pm 2\varepsilon_i: 1\le i \le n\}\,.
$$
We have an identification $Sp(V)/P=LG(V)$, the Lagrangian Grassmannian parametrizing all Lagrangian linear subspaces in $V$. 
It is an algebraic variety of dimension $n(n+1)/2$. 

We set $\rho(n)=(n,n-1,\ldots,2,1)$.

The Weyl group $W$ is here the hyperoctahedral group that can be identified with group of signed permutations, 
and $W_P=S_n$. The poset $W^P$ is naturally identified with the poset of strict partitions $\mu$ contained in $\rho(n)$ (see, e.g., \cite{Hi}),
and the corresponding Schubert class $Y^{\mu}$, where 
$$
\mu=(n\ge \mu_1> \cdots > \mu_r >0)
$$
is represented by the following variety. Consider a flag
$$
V_0\subset V_1 \subset \cdots \subset V_n \subset V
$$
of isotropic vector spaces with $\dim(V_i)=i$. Consider the following locus:
$$
\{L\in LG(V) : \dim(L\cap V_{n+1-\mu_i}) \ge i,  \ 1\le i \le r\}\,.
$$
The class of this locus does not depend on the flag, and is equal to the Schubert class $Y^\mu$.

\begin{theorem}\label{YQ} \cite{Pr1}, \cite{Pr2} \ In $A^{|\mu|}(LG(V))$, we have 
$$
Y^{\mu}=\\Q_{\mu} (R^{*}),
$$ 
where $R$ is the tautological subbundle on the Lagrangian Grassmannian.
\end{theorem}                  
The original argument \cite{Pr2} made use of a comparison of Pieri formulas for Lagrangian Schubert classes 
and for Schur $Q$-functions (for references, see \cite[Sect. 6]{Pr2}). There is also another proof in \cite[p. 40]{LaPr}, which uses the characteristic map 
for a Lagrangian Grassmannian, and relies on some divided difference and vertex operator computations.
\end{example}

We record the following ``duality'' result.
\begin{theorem} Let $G$ be a semisimple group, and let $P\subset G$ be a parabolic subgroup. For any $w\in W^P$ there exists 
exactly one $w'\in W^P$ such that $\dim X^w+\dim X^{w'}=\dim G/P$, and in $A^{*}(G/P)$ we have $X^w\cdot X^{w'}\ne 0$. In fact, 
$X^w\cdot X^{w'}=1$.
\end{theorem}
We call $X^{w'}$ the {\it dual class} to $X^w$. Let us discuss the following three examples.

For $P=B$, we have $w'=w_0\cdot w$. Indeed,
$$
X^w\cdot X^{w'}= \delta_{w,w_0 w'}X^{w_0}\,,
$$
where $X^{w_0}$ is the class of a point (see \cite[p. 20]{Ch}).

In the situation of Example \ref{Gr}, the dual class to $X^{\lambda}$ is  $X^{\lambda'}$ where $\l'_i=n-r-\lambda_{r+1-i}$ for $i=1,\ldots,r$. 
(See, e.g., \cite[p. 271]{Fu1}.) 

In the situation of Example \ref{LG}, the dual class to $Y^{\mu}$ is  $Y^{\mu'}$ where the parts of the strict partition $\mu'$ complement
the set of parts of $\mu$ in $\{1,2,\ldots,n\}$. (See \cite[p. 178]{Pr2}.) 

In general, the Poincar\'e duality for a partial flag manifold $G/P$ is described in terms of the Weyl group of $G$ in 
\cite[p. 197]{Kaji}.

\medskip
%{\small
The following Bertini-Kleiman theorem is often used to show positivity.
%\footnote{For Bertini's basic insight into this subject, see \cite{Ber}.}

\begin{theorem} \cite{Kl} \ Suppose a connected algebraic group $G$ acts transitively on a variety $X$ (over an algebraically
closed field $k$). Let $Y, Z$ be subvarieties of $X$. Then, denoting by $g\cdot Y$ the translate of $Y$ by $g\in G$, the following
two statements hold.

\noindent
(1) There exists a nonempty open subset $U\subset G$ such that for all $g\in U$, $(g\cdot Y)\cap Z$ is either empty
or of pure dimension \ 
$$
\dim (Y) +\dim (Z) - \dim (X)\,.
$$

\smallskip

\noindent
(2) If $Y$ and $Z$ are nonsingular, and char(k)=0, then there is a nonnempty open subset $U\subset G$ 
such that for all $g\in U$, $(g\cdot Y)\cap Z$ is nonsingular.
\end{theorem}

\begin{corollary} With the notation of the theorem, if $\dim(Y)+\dim(Z)=\dim(X)$, then $(g\cdot Y)\cap Z$ is either empty or 
a zero-dimensional subscheme. Under the assumptions of (2), all points in $(g\cdot Y)\cap Z$ are regular.
\end{corollary}

We end this section with the following fact on positivity of Schubert classes. Assertions (i) and (ii) are classically
known (also for $G/P$ with similar proofs).

\begin{proposition}\label{subsch} (i) Let $Z$ be a subvariety of $G/B$. If in $A^{*}(G/B)$, we have
$$
[Z]= \sum_{w\in W} a_w X^w\,,
$$
where $a_w\in \Z$, then all the coefficients $a_w$ are nonnegative.

\smallskip

\noindent
(ii) If for $w, v \in W$, in $A^{l(w)+l(v)}(G/B)$, we have
$$
X^w \cdot X^v=\sum_u c_{w v}^u X^u\,,
$$
then $c_{w v}^u \ge 0$.

\smallskip

\noindent
(iii)
Let $G\subset H$ be an inclusion of algebraic groups. Let $Q\subset H$
be a parabolic subgroup. Set $P=G\cap Q$, and let $i:G/P \to H/Q$ be the inclusion.
If $Z\subset H/Q$ is a subvariety, and in $A^{*}(G/P)$ we have
$$
i^*[Z]=\sum_{w\in W^P} a_w X^w\,,
$$
with $a_w \in \Z$, then all the coefficients $a_w$ are nonnegative.
\end{proposition}
Let us show, for instance, (i) and (iii). As for (i): For any $w$, we have 
$$
a_w=\int_{G/B} [Z]\cdot X^{w'}\,,
$$
where $X^{w'}$ is the dual class to $X^w$. Let $Y$ be a subvariety representing the class $X^{w'}$. We apply to $Z$ and $Y$ the Bertini-Kleiman theorem: for a general $g\in G$ we obtain a zero-dimensional scheme $(g\cdot Z)\cap Y$ and $a_w$ is its length, hence $a_w\ge 0$.

\smallskip

\noindent       
As for (iii):
We use the Bertini-Kleiman theorem for the subvarieties $Z$ and $G/P$ of $H/Q$: for a general $h\in H$, $h\cdot Z$ and $G/P$ meet
properly. Let 
$$
V= (h\cdot Z)\cap G/P\subset G/P\,,
$$ 
a schematic intersection. We now use Proposition \ref{subsch}(i) for the subvariety $V$ of $G/P$.
Alternatively, to conclude, we can use again the Bertini-Kleiman theorem, this time for $V\subset G/P$ and 
a subvariety representing the dual class to $X^w$.
\qed

\begin{corollary} Let $V$ be symplectic vector space of dimension $2n$ and let $LG(V)$ be the Lagrangian Grassmannian. 
Denote by
$$
i: LG(V) \hookrightarrow G_n(V)
$$ 
the inclusion. If in $A^{*}(LG(V))$ we have 
$$
i^*(X^\lambda)=\sum a_\mu Y^\mu\,,
$$ 
where $a_\mu \in \Z$, then $a_\mu \ge 0$.
\end{corollary}
A combinatorial positive rule for the coefficients $a_\mu$ was given in \cite[Prop. 2]{Pr3}.
%}
\section{Ample vector bundles and positive polynomials}\label{amppos}

In this section, we work over an algebraically closed field $k$ of arbitrary characteristic.
Let $X$  be a scheme, and let $\cF$ be a sheaf of $\cO_X$-modules. We say that $\cF$ is {\it generated by its global sections}
if there is a family of global sections $\{s_i\}_{i\in I}$, $s_i\in \Gamma(X,\cF)$, such that for each $x\in X$, the images of $s_i$
in the stalk $\cF_x$ generate that stalk as an $\cO_x$-module.

Recall that a line bundle $\cO(D)$ on a smooth curve of genus $g$ is generated by its global sections
if $\deg D\ge 2g$. It is ample iff $\deg D>0$; so a sufficiently high power of an ample line bundle is generated 
by its global sections. 

This is also the case of vector bundles of higher ranks on higher dimensional varieties.
Given a vector bundle $E$, we denote by $S^p(E)$ its $p$th symmetric power.
We say that a vector bundle $E$ on a variety $X$ is {\it ample} if for any sheaf $\cal F$ there exists $p_0\in \N$ such that  
for any $p\ge p_0$, the sheaf \ 
$
S^p(E)\otimes {\cal F}
$ 
is generated by its global sections. This is equivalent to say that the Grothendieck invertible sheaf ${\cal O}(1)$ on
$\P(E^*)$, the projective bundle\footnote{i.e. the bundle of lines in the fibers} of $E^*$ is ample.

Let us mention two properties of ample vector bundles \cite{Ha}: 
\begin{itemize}

\item a direct sum of ample vector bundles is ample;

\item for a partition $\lambda$, the {\it Schur bundle} $S^\lambda(E)$ (see \cite{La}, \cite[p. 131]{FuPr}) of an ample vector bundle $E$ is ample.
\end{itemize}

Perhaps this is a good moment to come back to positivity. Consider the following example.
If $E$ is a vector bundle, $\lambda$, $\mu$ partitions, then the integer coefficients $a_\nu$ in the expansion 
of the Schur polynomials of the Schur bundle $S^\lambda(E)$,
$$
s_\mu (S^\lambda(E))=\sum_\nu a_\nu s_\nu(E)\,,
$$ 
in the basis of Schur functions $\{s_\nu(E)\}$ are {\it nonnegative} \cite[Cor. 7.2]{Pr5} (see also \cite[Ex. 8.3.13]{Lz}).
This is a consequence of the second property. This information is nontrivial even for Chern classes (i.e. for $\mu=(1,\ldots,1)$);
for some examples of explicit computations, see \cite{Pr5}.

\smallskip

Let $E$ be a vector bundle of rank $n$ on a variety $X$, and $C$ a subscheme of $E$. We say that $C\subset E$ is a {\it cone} 
if it is stable under the natural $\G_m$ action on $E$.
If $C\subset E$ is a cone of pure dimension $c$, then one may intersect its cycle $[C]$ with the
zero-section of the vector bundle:
\begin{equation}\label{z}
z(C,E)=s_E^*([C])\in A_{c-n}(X)\,,
\end{equation}
where $s_E^*: A_c(E)\to A_{c-n}(X)$ is the Gysin map determined by the zero section $X\to E$.
In fact, we can use any other section $X\to E$ and $z(C,E)$ is the unique cycle class on $X$ such that
\begin{equation}\label{charz}
p^*(z(C,E)) = [C]
\end{equation}
in $A_c(E)$ (see \cite[(1.4)]{FuLa1}).

Here is an example of a positivity result with a pretty simple proof.

\begin{lemma} Let $E$ be a vector bundle on a variety $X$, and let $C$ be an irreducible cone in $E$.
If $E$ is generated by its global sections, then $z(C,E)$ is represented by a nonnegative cycle.
\end{lemma} 
\proof
Restricting $E$ to the support of $C$ \footnote{Cf. \cite[B.5.3]{Fu1}.}, we may assume that this support is equal to $X$. The inclusion $C\subset E$ 
gives rise to a subscheme $\P(C)\subset \P(E)$. If $E$ is generated by its global sections, then
$\cO(1)$ on $\P(E)$ is generated by its global sections. By the Bertini theorem, a general hypersurface
section on $\P(E)$ intersects $\P(C)$ properly or this intersection is empty. Hence a general section of $E$ intersects $C$
properly or the intersection is empty. Therefore $z(C,E)$ is represented by a nonnegative cycle.\qed 

For a projective variety $X$, there is well defined {\it degree} 
$$
\int_X: A_{0}(X)\to {\bf Z}
$$ 
(see \cite[Def. 1.4]{Fu1}). The following result of Fulton and Lazarsfeld is basic for applications to positivity.
  
\begin{theorem}\label{zCE}\cite{FuLa1} \ Let $E$ be an ample vector bundle of rank $n$ on a projective variety $X$. Let $C$ be a cone in $E$ 
of pure dimension $n$. Then we have 
$$
\int_X z(C,E)>0\,.
$$
\end{theorem}
For a more extensive study of positivity in intersection theory, coming from ample vector bundles and vector bundles generated by their
global sections, see \cite[Thm 12.1]{Fu1}. 

\begin{remark}\label{zh} \rm Suppose $k=\C$. Under the assumptions of Theorem \ref{zCE}, we have in $H_0(X,{\bf Z})$ 
the homology analog of $z(C,E)$, denoted by the same symbol.
We also have the homology degree map $\deg_X: H_0(X,{\bf Z})\to {\bf Z}$.
They are compatible with their Chow group counterparts via the cycle map: $A_0(X) \to H_0(X,{\bf Z})$ 
(cf. \cite[Sect. 19]{Fu1}). Thus we have
\begin{equation}\label{ineq}
\deg_X \bigl(z(C,E)\bigr)>0\,.
\end{equation}
\end{remark}

We record the following result.
\begin{proposition}\label{CV} \cite{MPW1}
Let $E$ be a vector bundle on a complete homogeneous variety $X$. Let $C$ be a cone in $E$ and let $Y\subset X$
be a subvariety of dimension $\dim(X)+\rank(E) -\dim C$. Assume that $E$ is generated by its global sections. 
Then the intersection \ $[C]\cdot [Y]$ is nonnegative.
\end{proposition}

\medskip

Let $c_1,c_2,\ldots$ be commuting variables with $\deg(c_i)=i$. Fix $d, n \in {\bf N}$. Let $P(c_1,\ldots,c_n)$ be a weighted 
homogeneous polynomial of degree $d$. We say that $P$ is {\it numerically positive for ample vector bundles}, or simply {\it positive}, 
if for every $d$-dimensional 
projective variety $X$ and any ample vector bundle of rank $n$ on $X$, we have 
$$
\int_X P(c_1(E),\ldots, c_n(E)) > 0\,.
$$

\smallskip

For example, Griffiths \cite{Gr} who pioneered this subject, found the following positive polynomials: $c_1$, $c_2$, $c_1^2-c_2$. 
Bloch-Gieseker \cite{BlGi} showed that $c_d$ is positive for $d\le n$.

Given a partition $\lambda$, with the conjugate partition $\mu$, we set
\begin{equation}\label{sl}
s_\lambda=s_\lambda(c_1, c_2,\ldots):=|c_{\mu_i-i+j}|_{1\le i,j \le l(\mu)}\,.
\end{equation}
Kleiman \cite{Kl1} showed that positive polynomials for surfaces are nonnegative combinations of $s_2$ and $s_{1,1}$.
Gieseker \cite{Gi} proved that $s_d$ (the $d$-th Segre class) is positive.

Fulton and Lazarsfeld gave the following characterization of positive polynomials. Let $P$ be a weighted homogeneous polynomial 
of degree $d$ in $n$ variables. Write
\begin{equation}\label{al}
P=\sum_\lambda a_\lambda s_{\lambda}\,,
\end{equation}
where $a_\lambda \in \Z$.
       
\begin{theorem}\label{pos} \cite{FuLa1} The polynomial $P$ is positive iff $P$ is not zero and all the coefficients $a_\lambda$
in (\ref{al}) are nonnegative. 
\end{theorem}

The proof of the theorem combines the hard Lefschetz theorem appropriately adapted to this subject by Bloch and Gieseker \cite{BlGi} 
and the Giambelli formula, which was recalled in Theorem \ref{Gia}. 

\begin{remark} \rm We now discuss some results related to Theorems \ref{zCE} and  \ref{pos}. The latter was generalized by Demailly, Peternell
and Schneider to {\it nef} vector bundles in \cite{DPS}. The former has a very simple proof due to Fulton and Lazarsfeld in \cite{FuLa2} 
in the case when $E$ is ample and generated by its global sections. Hacon \cite{Hac} showed that 
these assumptions are not sufficient, for a positive polynomial $P$, to have 
$
\int_X P(E)\ge P(n, {n \choose 2},\ldots,{n \choose d})\,,
$ 
as it was conjectured by Beltrametti, Schneider and Sommese in \cite{BSS}. This last inequality is true for {\it very ample} bundles 
(${\it loc.cit.}$). 
Consider a vector bundle $E$ on a complex projective manifold. Griffiths \cite{Gr} defined $E$ to be {\it numerically positive} if for any analytic subvariety $W\subset M$, and any rank $q$ quotient $Q$ of $E_{|W}$, we have $\int_W P(c(Q)) >0$ for any homogeneous polynomial of degree equal to $\dim(W)$ from the Griffiths cone associated to $q$ (see also \cite[App. A]{FuLa1}).
Griffiths speculated on the possibility that arbitrary ample bundles are numerically positive. This was proved, using Schubert calculus, by Usui and Tango \cite{UsTa} for bundles generated by their global sections. The numerical positivity of all ample bundles was proved in \cite[App. A]{FuLa1}.
\end{remark}

\section{Thom polynomials for singularities of maps}\label{tsm}

Thom polynomials came from algebraic topology and singularities. They are tools to measure the complexity of singularities. 
In this section, we investigate Thom polynomials of singularities of maps. Let 
$$
f:M\to N
$$ 
be a map of complex analytic manifolds; 
we say that $x\in M$ is a singularity of $f$ if $df_x$ fails to have the maximal rank. 

We now follow the terminology from \cite{Rim} for what concerns map germs $({\bf C}^m,0) \to ({\bf C}^n,0)$ and their stable versions.
Two map germs $\kappa_1, \kappa_2 :({\bf C}^{m},0)\to ({\bf C}^{n},0)$ are said to be {\it right-left equivalent} if there exist germs of biholomorphisms $\phi$ of $({\bf C}^m,0)$ and  $\psi$ of $({\bf C}^{n},0)$ such that $\psi \circ \kappa_1 \circ \phi^{-1} = \kappa_2.$  A {\it suspension} of a germ map $\kappa$ is its trivial unfolding: $(x,v)\mapsto (\kappa(x),v)$. Let us fix $l\in {\bf N}$.  
Consider the equivalence relation on stable map germs $({\bf C}^{\bullet},0)\to ({\bf C}^{\bullet+l},0)$ generated by right-left equivalence and suspension. An equivalence class of this relation is
often called {\it singularity} and denoted by $\eta$. 

According to Mather's classification \cite{dPW}, the finite dimensional (local) ${\bf C}$-algebras are in one-to-one correspondence with classes of 
contact equivalence classes of singularities (cf. \cite{FeRi}). For instance,
$A_i$ stands for the stable germs with local algebra ${\bf C}[[x]]/(x^{i+1})$, $i\ge 0$;
and $I_{a,b}$ for stable germs with local algebra 
${\bf C}[[x,y]]/(xy, x^a+y^b)$, \ $b\ge a\ge 2$ (they also depend on $l$).

Following Thom, we look at the locus
$$
\eta(f):=\{x\in M : \hbox{the singularity of} \ f \ \hbox{at} \ x \
\hbox{is} \ \eta \}
$$
for a map $f: M\to N$, and try to compute its class in terms of the Chern classes of $M$ and $N$.
More precisely, we restrict ourselves only to {\it general} maps, i.e., the maps from some open subset in
the space of all maps.

For example, let $f:M\to N$ be a general morphism of compact Riemann surfaces. Suppose that the singularity 
is $A_1$: $z\mapsto z^2$. Then $\eta(f)$ is the ramification divisor of $f$, and by the Riemann-Hurwitz formula the wanted
class is $f^ *c_1(N)-c_1(M)$. We refer the reader to \cite[p. 300]{Kl3} for more details.

The space of germs of maps from $(\C^m,0)$ to $(\C^n,0)$ has infinite dimension, which is inconvenient from the point 
of view of algebraic geometry. To remedy this, we pass to the {\it spaces of jets} of germs of maps.
Fix $m,n,p \in {\bf N}$. Consider the space 
$
{\cJ}^p({\bf C}^m_0,{\bf C}^n_0)
$ 
of $p$-jets of analytic functions from ${\bf C}^m$ to ${\bf C}^n$, which map $0$ to $0$ (see \cite[pp. 36-38]{AGV}). This space will also be denoted by $\cJ(m,n)$ or simply by $\cJ$ to avoid too complicated notation.

Let $\Aut_n^p$ denote the group of $p$-jets of automorphisms of $({\bf C}^n,0)$.  

Consider the natural right-left action of the group $\Aut_m^p \times \Aut_n^p$ on the space ${\cJ}^p({\bf C}^m_0,{\bf C}^n_0)$. 
By a {\it singularity class} we mean a closed algebraic right-left invariant subset of ${\cJ}^p({\bf C}^m_0,{\bf C}^n_0)$. 

Given complex analytic manifolds $M^m$ and $N^n$, a singularity class $\Sigma \subset {\cJ}^p({\bf C}^m_0,{\bf C}^n_0)$ defines 
the following subset $\Sigma(M,N)\subset {\cJ}^p(M,N)$, where ${\cJ}^p(M,N)$ is the space of $p$-jets of maps from $M$ to $N$:
using the coordinate maps $M\cong{\bf C}^m$ and $N\cong{\bf C}^n$, we declare that a point belongs to $\Sigma(M,N)$ iff it belongs 
to $\Sigma$. If we change the coordinate maps, then the set $\Sigma(M,N)$ remains unchanged by virtue of right-left invariance of $\Sigma$.

\begin{theorem}\label{Thompol}
Let $\Sigma \subset {\cJ}^p({\bf C}^m_0,{\bf C}^n_0)$ be a singularity class.
There exists a universal polynomial ${\cT}^{\Sigma}$ over $\bf Z$
in $m+n$ variables $c_1,\ldots,c_m,c'_1,\ldots,c'_n$ which depends only on $\Sigma$, $m$ and $n$
such that for any manifolds $M^m$, $N^n$ and for a general map $f: M\to N$,
the class of 
$$
\Sigma(f):=(j^pf)^{-1}(\Sigma(M,N))
$$ 
is equal to
$$
{\cal T}^{\Sigma}(c_1(M),\ldots,c_m(M),f^*c_1(N),\ldots,f^*c_n(N)),
$$
where $j^pf: M\to \cJ^p(M,N)$ is the $p$-jet extension of $f$ (\cite[pp. 36-38]{AGV}).
\end{theorem}
This is a theorem due to Thom, see \cite{Th}. The polynomial $\cT^\Sigma$ is called the {\it Thom polynomial of} $\Sigma$.

Note that a singularity $\eta$ corresponds here to the singularity class $\Sigma$ being the closure of a single right-left
orbit, and the locus $\eta(f)$ is generalized to $\Sigma(M,N)$. The key problem is to compute the classes of
these varieties in terms of the Chern classes of the manifolds $M$ and $N$.

\begin{lemma}\label{scc} Let $\Sigma\subset \cJ$ be a singularity class. Then $\Sigma$ is a cone in the vector space $\cJ$.
\end{lemma}
\proof For a function $f\in \Sigma$ and a scalar $c\in \C^*$, we have $c\cdot f\in \Sigma$ because $\G_m\subset \Aut_n^p$, and the
singularity class $\Sigma$ is $\Aut_m^p\times \Aut_n^p$-invariant.\qed

\medskip

We follow Kazarian's approach to Thom polynomials of singularities of maps \cite{Ka1}. We set
$$
G:={\Aut}_m^p \times {\Aut}_n^p \,.
$$
Consider the classifying principal $G$-bundle $EG\to BG$  \ \cite{Mi} (see also \cite[Sect. 7]{H}). 
Here $EG$ is a contractible space with a free action of the group $G$. This action extends to the diagonal
action on the product space $EG \times \cJ$. 
Invoking \cite[Def. 3.1]{H} and its notation, we set
$$
\widetilde{\cJ}:=EG\times^G \cJ=(EG \times \cJ)/G\,.
$$
This space is often called the {\it classifying space of singularities}\footnote{Note that the same construction
is used in the definition of Borel of equivariant cohomology for a $G$-space $\cJ$.}. 
For a given singularity class $\Sigma \subset \cJ$, we define
$$
\widetilde{\Sigma}:=EG\times^G \Sigma\subset\widetilde{\cJ}\,.
$$
We have $\codim(\widetilde{\Sigma},\widetilde{\cJ})= \codim(\Sigma,{\cJ})$. We denote by ${\cT}^{\Sigma}\in H^{2\codim(\Sigma,\cJ)}(\widetilde{\cJ},{\bf Z})$ the dual class of $[\widetilde{\Sigma}]$. 
The classifying spaces $BG$, $\widetilde{\cJ}$, etc. have infinite dimensions and the notion of the ``dual class''
should be clarified, see \cite[Rem. 1.6]{Ka1} and \cite[footnote $(^{4})$]{PrWe1}.
%\footnote{ 
%We may approximate $EG\to BG$ by a sequence of $G$-bundles over finite dimensional manifolds $(EG)_i\to (BG)_i$, where 
%$i\to \infty$. Then ``${\cT}^{\Sigma}$ is the dual class to the fundamental class of $\widetilde{\Sigma}$''
%means that for any $i$ the image of ${\cT}^{\Sigma}$ in $H^{*}((BG)_i, {\bf Z})\cong H^{*}((EG)_i\times^G\cJ,{\bf Z})$ 
%is Poincar\'e dual to the fundamental class of $(EG)_i\times^G\Sigma$ in $H_{*}((EG)_i\times^G\cJ, {\bf Z})$, where
%the last homology group is the Borel-Moore homology (see, e.g., \cite[Sect. 19]{Fu1}).}. 

The projection to the second factor $\widetilde{\cJ}\to BG$ is a bundle with fiber isomorphic to $\cJ$ and structure group $G$.
Since $\cJ$ is contractible, and also $G$ is contractible to the subgroup $GL_m\times GL_n$ of linear changes, we get
$$
H^*(\widetilde{\cJ}, {\bf Z})
\cong H^*(BG,{\bf Z})\cong H^*(BGL_m\times BGL_n,{\bf Z})\,.
$$
Hence ${\cT}^{\Sigma}$ is identified with a polynomial in $c_1,\ldots, c_m$
and $c_1',\ldots, c_n'$ which are the Chern classes of universal bundles
$R_m$ and $R_n$ on $BGL_m$ and $BGL_n$.
This is the Thom polynomial $\cT^\Sigma$.

We now explain what we mean by {\it stable} singularity class. The suspension
$$
\cS:\cJ(m,n)\hookrightarrow \cJ(m+1,n+1)
$$
allows one to increase the dimension of the source and the target
simultaneously: with the local coordinates $x_1,x_2,\ldots$ for the source
and a function $f=f(x_1,\ldots,x_m)$, the jet
${\cS}(f)\in\cJ(m+1,n+1)$ is defined by
$$
{\cS}(f)(x_1,\dots,x_m,x_{m+1}):=(f(x_1,\dots,x_m),x_{m+1})\,.
$$
Suppose that the singularity class $\Sigma$ is
{\it stable under suspension}. By this we mean that it is a member
$\Sigma_0=\Sigma$ of a family
$$
\{\Sigma_r\subset\cJ(m+r,n+r)\}_{r\ge 0}
$$
such that
$$
\Sigma_{r+1}\cap\cJ(m+r,n+r)=\Sigma_r
$$
and
$$
{{\cT}^{\Sigma_{r+1}}}_{|H^*(BGL_{m+r}\times BGL_{n+r},{\bf Z})}
={\cT}^{\Sigma_r}\,.
$$
This means that if we specialize
$$
c_{m+r+1}=c'_{n+r+1}=0
$$
in the polynomial ${\cT}^{\Sigma_{r+1}}$, we obtain the polynomial
${\cT}^{\Sigma_r}$. If $\Sigma$ is closed under the {\it contact equivalence} (see \cite{FeRi}), then it is
stable in our sense.

The theorem of Thom has the following refinement due to Damon
\cite{Da} for singularity classes $\Sigma$ which are stable
under suspension: ${\cT}^{\Sigma}$ is a polynomial in
$$
c_i(R_m-R_n) \,, \ \ \ \hbox{where} \ \ \ i=1,2,\ldots\,. 
$$
So, we can use the bases of monomials in the Chern classes in $R_m-R_n$ or $R_n-R_m$ or $R_m^*-R_n^*$ or $R_n^*-R_m^*$.
We can also use the bases of (supersymmetric) Schur functions in $R_m-R_n$ or $R_n-R_m$ or $R_m^*-R_n^*$ or $R_n^*-R_m^*$.
About a decade ago, calculations of the Thom polynomials using the bases of Schur functions were done independently by 
Feh\'er-K\"om\"uves and Lascoux-Pragacz. 

For Morin singularities $A_i$, there is a positivity conjecture of Rim\'anyi (1998), asserting that the expansions of the 
Thom polynomials $\cT^{A_i}$ in the basis of monomials in the Chern classes in $R_n-R_m$ have nonnegative coefficients. 
See \cite{Be} for a discussion of a link of this conjecture with the Green-Griffiths conjecture about holomorphic curves in
nonsingular projective varieties.

\begin{example} \rm We display here three Thom polynomials for the Morin singularities between equal dimensional manifolds (so $l=0$ in
the notation from the beginning of this section):

\noindent
$A_3$: \  $c_1^3+3c_1c_2+2c_3$;

\noindent
$A_4$: \ $c_1^4+6c_1^2c_2+2c_2^2+9c_1c_3+6c_4$;

\noindent
$A_5$: \ $c_1^5+10c_1^3c_2+25c_1^2c_3+10c_1c_2^2+38c_1c_4+12c_2c_3+24c_5$.

\end{example}

In general, the expansions of Thom polynomials of stable singularities in the basis of monomials in the Chern classes of $R_n-R_m$ 
can have negative coefficients.

\begin{example} \rm We give here three Thom polynomials for the singularities $I_{p,q}$ between equal dimensional manifolds
(for the definition of these singularities, see the beginning of this section):

\noindent
$I_{2,2}$: \ $c_2^2-c_1c_3$;

\noindent
$I_{2,3}$: \ $2c_1c_2^2-2c_1^2c_3+2c_2c_3-2c_1c_4$;

\noindent
$I_{2,4}$: \ $2c_1^2c_2^2+3c_2^3-2c_1^3c_3+2c_1c_2c_3-3c_3^2-5c_1^2c_4+9c_2c_4-6c_1c_5$.

\end{example}

It is not obvious that $\cT^{\Sigma}\ne 0$ for a nonempty stable singularity class $\Sigma$.

We now examine the expansions of Thom polynomials of stable singularities in the basis $\{s_\lambda(R_n-R_m)\}$ labeled by partitions $\lambda$.
We refer the reader for a variety of examples to \cite[p. 93-94]{PrWe1}, \cite{FeRi} and \cite{OzPr}.\footnote{In \cite{PrWe1} 
and \cite{OzPr} the authors worked with the basis of Schur functions $\{s_{\lambda}(R_m^*-R_n^*)\}$, so the Schur functions given 
in the examples there are labeled by the conjugate partitions of those appearing in the present convention.}

\medskip

\begin{theorem}\label{stab} \cite{PrWe1} \ Let $\Sigma$ be a nonempty stable singularity class. Then
for any partition $\lambda$ the coefficient $a_\lambda$ in 
\begin{equation}\label{TMTN}
{\cal T}^{\Sigma}= \sum a_\lambda s_\lambda(R_n-R_m) 
\end{equation}
is nonnegative and $\sum a_\lambda>0$.
\end{theorem}
This feature of Schur function expansions of Thom polynomials was pointed out in \cite{Pr6}, 
conjectured for Thom-Boardman singularities by Feh\'er and K\"om\"uves \cite{FeKo} (they computed the 
Schur function expansions of the Thom polynomials of $\Sigma^{i,j}[-i+1]$),
and conjectured for all singularity classes in \cite{Pr7}. 

\smallskip

Note that each partition $\l$ appearing in the RHS of (\ref{TMTN}) is contained in the $(n,m)$-hook (see \cite[p. 35]{FuPr}).

To prove the theorem, we generalize the equation (\ref{TMTN}) for any pair of complex vector bundles ($E$, $F$) on any 
paracompact space $X$. 
To this end, we apply the techniques of fiber bundles. Apart from vector bundles, we also use principal $G$-bundles 
associated with finite collections of vector bundles\footnote{The associated principal $GL_n$-bundle of a vector bundle $E$ of rank $n$ is often called the {\it frame bundle} of $E$ (its fibers consist of all ordered bases of the fibers of $E$).} on a common base space (here $G=\prod_i GL_{n_i}$, where $n_i$ are the ranks of the vector bundles). For principal bundles, we refer, e.g., to \cite[Sect. I.5]{KoNo} 
or \cite[Sect. 5]{H}.

Moreover, it is convenient to pass to the topological homotopy category, where any pair of vector bundles can be pulled back
from the universal pair of vector bundles on $BGL_m\times BGL_n$.

We first pull back the bundle $\widetilde{\cJ}$ from $BG$ to $BGL_m\times BGL_n$ via 
the embedding
$$
GL_m \times GL_n \hookrightarrow {\Aut}_m\times {\Aut}_n\,.
$$
Since $GL_m \times GL_n$ acts linearly on $\cJ$, the obtained
pullback bundle is now the following vector bundle on $BGL_m\times BGL_n$:
$$
\cJ(R_m,R_n):=\Bigl(\bigoplus_{i=1}^p S^i(R_m^*)\Bigr) \otimes R_n\,.
$$
The bundle $\cJ(R_m,R_n)$ contains the preimage of
$\widetilde{\Sigma}$, denoted by $\Sigma(R_m,R_n)$, whose class is
\begin{equation}\label{srr}
[\Sigma(R_m,R_n)]=\sum_\lambda a_\lambda s_\lambda(R_n-R_m)\,,
\end{equation}
with the same coefficients $a_\lambda$ as in (\ref{TMTN}).

Consider now a pair of vector bundles $E$ and $F$ of ranks $m$ and $n$ on a variety $X$.
We set
$$
\cJ(E,F):=\Bigl(\bigoplus_{i=1}^p S^i(E^*)\Bigr) \otimes F\,.
$$

Let $P(E,F)$ be the principal $GL_m\times GL_n$-bundle
associated with the pair of vector bundles $(E,F)$. We have
$$
\cJ(E,F)=P(E,F)\times^{GL_m\times GL_n}\cJ\,.
$$
We set
$$
\Sigma(E,F):=P(E,F)\times^{GL_m\times GL_n}\Sigma\subset\cJ(E,F)\,,
$$
a locally trivial fibration with the fiber equal to $\Sigma$. 

\begin{lemma}\label{cEF} The variety $\Sigma(E,F)$ is a cone in the vector bundle $\cJ(E,F)$.
\end{lemma}
\proof The assertion follows from Lemma \ref{scc}.\qed

\smallskip

\begin{lemma}\label{lEF} The dual class of $[\Sigma(E,F)]\in H_{2\dim(\Sigma)}(\cJ(E,F),{\bf Z})$ in
$$
H^{2\codim(\Sigma,\cJ)}(\cJ(E,F),{\bf Z})=H^{2\codim(\Sigma,\cJ)}(X,{\bf Z})
$$
is equal to
\begin{equation}\label{EF}
\sum_\lambda a_\lambda s_\lambda(F-E)\,,
\end{equation}
where the coefficients $a_\lambda$ are the same as in (\ref{TMTN}) \footnote{The meaning of the ``dual class of $[\Sigma(E,F)]$'' 
for a {\it singular} $X$ is explained in \cite[Note 6]{PrWe1}.}.
\end{lemma} 
\proof The pair of vector bundles ($E$,$F$) on a variety $X$ can be pulled back from the universal pair ($R_m$, $R_n$) 
on $BGL_m \times BGL_n$ using a $C^\infty$ map. We get the assertion of the lemma by pulling back the equation (\ref{srr}).
Consequently, the coefficients of $s_\lambda(F-E)$ in (\ref{EF}) are the same as the coefficients of $s_\lambda(R_n-R_m)$
in (\ref{srr}).\qed 

\bigskip

\noindent
{\bf Proof of Theorem \ref{stab}.}\footnote{This is the same proof as that in \cite{PrWe1}, but with ``mehr Licht''.}
Let $e=\codim(C,\cJ)$. This means that for any partition $\lambda$ appearing in (\ref{EF}) its weight $|\l|$ is equal to $e$. 

The idea of the proof is to produce from (\ref{EF}) a numerically positive polynomial for ample vector bundles, which captures
positivity information about all the $a_\lambda$'s. Since (\ref{EF}) is a supersymmetric polynomial, and we want
a usual symmetric polynomial, we wish to specialize $E$ to be a trivial bundle. Since the singularity class $\Sigma$ is stable, 
we can use a pair of vector bundles $E$ and $F$ on $X$ of the corresponding ranks $m'=m+r$ and $n'=n+r$ for some $r \ge 0$, 
instead of $m$ and $n$. So we can assume that $n'>>0$. In particular, we may suppose that $n'\ge e$. 

We use a specialization argument: let $X$ vary over projective varieties of dimension $e$, let $F$ vary over ample vector
bundles of rank $n'$ on $X$, and let $E$ be a trivial vector bundle ${\bf 1}^{m'}$ of rank $m'$ on $X$. By the theory of symmetric functions 
(see, e.g., \cite[Sect. 3.2]{FuPr}), the Schur polynomials $s_{\l}(F-E)$ appearing in (\ref{EF}) are indexed by partitions 
$\l$ of weight $|\l|=e$, which are contained in the $(n,m)$-hook. In general, such polynomials vanish under our specialization. 
But the assumption $n'\ge e$, or equivalently, $\rank F\ge |\lambda|$, guarantees that the partition corresponding to a summand 
$a_\lambda s_{\l}(F-E)$ appearing in (\ref{EF}) has at most $n'$ parts, and thus this summand survives the specialization, 
giving $a_\lambda s_{\l}(F)$. After the specialization, the expression (\ref{EF}) becomes
\begin{equation}\label{F}
\sum_\lambda a_\lambda s_\lambda(F)\,,
\end{equation}
where the summation is as in (\ref{EF}). Consider the polynomial 
$$
P:= \sum_\lambda a_\lambda s_\lambda\,,           
$$
with the $s_\lambda$'s as in (\ref{sl}) and the summation as in (\ref{F}). We want to show that $P$ is positive.
To this end, consider the cone $\Sigma(E,F)$ in $\cJ(E,F)$ (see Lemma \ref{cEF}) and its cone class 
$z(\Sigma(E,F),\cJ(E,F))$ (see (\ref{z}) and Remark \ref{zh}). Since  $\dim \Sigma(E,F)= \rank \cJ(E,F)$, this cone class belongs to $H_0(X,\Z)$.
It follows from Lemma \ref{lEF} that the dual class of
$z(\Sigma({\bf 1}^{m'},F),\cJ({\bf 1}^{m'},F))$ is the element of $H^{2e}(X,\Z)$ given by the expression (\ref{F}).

Since a direct sum of ample vector bundles is ample (see \cite[Prop. (2.2)]{Ha}), and the 
vector bundle $\cJ({\bf 1}^{m'},F)$ is a direct sum of several copies of $F$, then $\cJ({\bf 1}^{m'},F)$ is ample.
Therefore by Theorem \ref{zCE} and the inequality (\ref{ineq}), we have
$$
\int_X P(F)=\deg_X(z(\Sigma({\bf 1}^{m'},F),\cJ({\bf 1}^{m'},F))>0\,,
$$
and thus conclude that $P$ is positive.

In turn, by Theorem \ref{pos} we get that $P$ is nonzero, and all the coefficients $a_{\l}$ are nonnegative; 
hence also $\sum_\lambda a_\lambda >0$.\qed

\medskip

\noindent
{\bf Question.} Does there exists a basis different (up to rescaling) from the basis $\{s_\lambda(R_n-R_m)\}$ with the property that any Thom
polynomial of a stable singularity class has a positive expansion in that basis?

%\end{document}     

\section{Thom polynomials for invariant cones}\label{tic}

In the previous section, in the context of classical Thom polynomials, we have investigated
the functor of $p$-jets:
$$
(E,F)\mapsto {\cJ}^p(E,F)=\Bigl(\bigoplus_{i=1}^p S^i(E^*)\Bigr)\otimes F\,,
$$
defined on pairs of vector bundles, where $p$ is large enough.

We want to generalize this setting.
Suppose that $(n_1,\ldots,n_l)\in {{\bf N}^*}^l$ and that $V$ is a representation of
$
G=\prod_{i=1}^l GL_{n_i}\,.
$
The representation $V$ gives rise to a {\it functor} $\phi$ defined for
a collection of bundles on a variety $X$:
$$
E_1,\dots,E_l\mapsto \phi (E_1,\ldots, E_l)\,,
$$
with $\dim E_i=n_i$, $i=1,\ldots, l$. By passing to the dual bundles,
we may assume that the functor $\phi$ is covariant in each variable.

Let $P(E_{\bullet})=P(E_1,\ldots,E_l)$ be the principal $G$-bundle associated with the vector bundles 
$E_1,\ldots,E_l$. We define a new vector bundle:
$$
V(E_{\bullet})=V(E_1,\ldots,E_l):=P(E_{\bullet})\times^G V
$$
with fiber equal to $V$.

Suppose now that a $G$-invariant cone $\Sigma\subset V$ is given.
We set
$$
\Sigma(E_{\bullet})=\Sigma(E_1,\ldots,E_l)
:=P(E_{\bullet})\times^{G}\Sigma\subset V(E_{\bullet})\,,
$$
a fibration with fiber equal to $\Sigma$.

Let $R^{(i)}$, $i=1,\dots,l$, be the pullback of the tautological vector
bundle from $BGL_{n_i}$ to 
$$
BG=\prod_{i=1}^l BGL_{n_i}\,.
$$

We denote by 
$$
\cT^\Sigma\in H^{2codim(\Sigma,V)}\bigl(V(R^{(1)},\ldots, R^{(l)}),{\bf Z}\bigr)=H^{2codim(\Sigma,V)}(BG,{\bf Z})
$$
the dual class\footnote{Here the ``dual class'' has the same meaning as in the approach to Thom polynomials via classyfying
spaces in the previous section.} of $[\Sigma(R^{(1)},\ldots, R^{(l)})]$, and call it the {\it Thom polynomial} of $\Sigma$.

Then, the so defined Thom polynomial ${\cT}^\Sigma\in H^{*}(BG, {\bf Z})$
depends on the Chern classes of the universal bundles $R^{(i)}$'s.
We write \ ${\cT}^{\Sigma}(E_1,\ldots,E_l)$ \ for
the Thom polynomial ${\cT}^{\Sigma}$, with $c_j(R^{(i)})$ replaced by
$c_j(E_i)$ for $i=1,\ldots,l$.

Arguing like in Lemma \ref{lEF}, we know that for any vector bundles $E_1,\ldots,E_l$
on a variety $X$, the class $[\Sigma(E_{\bullet})]$ in
$$
H^{2\codim(\Sigma,\cJ)}(V(E_{\bullet}), {\bf Z})=H^{2\codim(\Sigma,\cJ)}(X,{\bf Z})
$$
is equal to ${\cT}^{\Sigma}(E_1,\ldots,E_l)$.

Since the Schur functions form an additive basis of the ring of symmetric
functions, the Thom polynomial ${\cT}^\Sigma$ is uniquely written
in the following form:
\begin{equation}\label{bul}
{\cT}^\Sigma=\sum a_{\lambda^{(1)},\dots,\lambda^{(l)}} \ 
s_{\lambda^{(1)}}(R^{(1)})\cdot \ldots \cdot s_{\lambda^{(l)}}(R^{(l)})\,,
\end{equation}
where $a_{\lambda^{(1)},\dots,\lambda^{(l)}}\in \Z$ 

\smallskip

We say that the functor $\phi$, associated with a representation $V$,
{\it preserves spannedness} if for a collection of 
vector bundles $E_1,\dots,E_l$ generated by their global sections, the bundle
$
\phi(E_1,\dots,E_l)
$
is generated by its global sections.

Examples of functors preserving spannedness over fields of
characteristic zero are {\it polynomial functors}. They are, at the same time,
quotient functors and subfunctors of the tensor power functors (cf. \cite{Ha}).

\begin{theorem}\label{Tm} \cite{PrWe2} Suppose that the functor $\phi$
preserves spannedness. Then the coefficients $a_{\lambda_1,\ldots,\lambda_l}$ in
(\ref{bul}) are nonnegative. Assume additionally
that there exists a projective variety $X$ of dimension greater than or equal
to $\codim(\Sigma,V)$, and there exist vector bundles $E_1,\ldots,E_l$ on $X$
such that the bundle $\phi(E_1,\dots,E_l)$ is ample. Then at least one
of the coefficients $a_{\lambda_1,\ldots,\lambda_l}$ is positive.
\end{theorem}

\medskip

Consider now the Thom polynomial $\cT^\Sigma$ associated with a nonempty, possibly nonstable singularity class $\Sigma$
in the space of jets $\cJ(m,n)$.
By the theory of symmetric functions (see, e.g., \cite[Sect. 3]{FuPr}), there exist coefficients $b_{\lambda \mu}\in {\bf Z}$ such that
\begin{equation}\label{beta}
{\cT}^{\Sigma}=\sum_{\lambda,\mu} b_{\lambda \mu} s_\lambda(R_n)\cdot s_\mu(R_m^*)\,.
\end{equation}

The following result follows from Theorem \ref{Tm}.

\begin{corollary}\label{CpS} For any pair of partitions $\lambda,\mu$, we have
$b_{\lambda \mu}\ge 0$ and $\sum b_{\lambda \mu}>0$.
\end{corollary}

Let now $\Sigma$ be a stable singularity class.
There exist coefficients $a_\lambda \in {\bf Z}$ such that
\begin{equation}\label{alpha}
{\cT}^{\Sigma}=\sum_\lambda a_\lambda s_\lambda (R_n-R_m)\,,
\end{equation}
the sum is over partitions $\lambda$ with $|\lambda|=\codim(\Sigma,\cJ(m,n))$.

Here is another proof of Theorem \ref{stab}.
By the theory of symmetric functions ({\it loc.cit.}), we have that the coefficient of
$s_{\l}(R_n-R_m)$ in the RHS of (\ref{alpha}) is equal to the coefficient
of $s_{\l}(R_n)$ in the RHS of (\ref{beta}), that is, \ $a_\lambda=b_{\lambda,\emptyset}$ \ for
any partition $\l$. The assertion now follows from Corollary \ref{CpS}.

\begin{remark} \rm Another proof of the {\it nonnegativity} assertions of Theorem \ref{stab} and Corollary \ref{CpS} 
was communicated to the author by Klyachko and indendependently by Kazarian (private communications). 
For details, see \cite[p. 452]{OzPr}. These proofs use the Bertini-Kleiman theorem.  
Coming back to the above proofs of Theorem \ref{stab}, we see that the use of ample vector bundles and 
the Fulton-Lazarsfeld Theorem \ref{zCE}, apart from the nonnegativity of the considered coefficients, implies 
that at least one of them is strictly positive.

\end{remark}

\section{Lagrangian Thom polynomials}\label{tlag}

Lagrangian Thom polynomials were considered by Vassiliev \cite{Va} (see also \cite{Ka2}).

Let us fix a positive integer $n$. Suppose that $W$ be a complex vector space, where $\dim W=n$. Let
$$
V=W\oplus W^*
$$
be a linear symplectic space, equipped with the symplectic form $\langle , \rangle$, defined by
$$
\langle (w_1,f_1),(w_2,f_2)\rangle=f_1(w_2)-f_2(w_1)
$$
for $w_i \in W$ and $f_i \in W^*$, $i=1,2$. We view $V$ as a symplectic manifold. Writing $q=(q_1,\ldots,q_n)$ 
for the coordinates of $W$ and $p=(p_1,\ldots,p_n)$ for the dual coordinates of $W^*$, the symplectic form 
on $V$ is \ $\sum_{i=1}^n dp_i \wedge dq_i$. 

Denote by \ $\varrho: V \to W$ the projection.

Any germ of a Lagrangian submanifold $L$ of $V$ through $0$ such that $\varrho_{|_L}$ is a submersion
is a graph of a 1-form \ $\alpha: W\to W^*$.
The condition that $L$ is Lagrangian is equivalent to
  \ $d\alpha=0$.
Since we deal with germs, we can write
  \ $\alpha=df$ \
for some function $f:W \to {\C}$.

The space of germs of Lagrangian submanifolds $L\subset V$ passing through $0$ 
has infinite dimension, which is inconvenient from the point of view of algebraic geometry. 
To remedy this, we pass to the space of jets of germs of Lagrangian submanifolds.

Let us fix, once for all, a nonnegative integer $p$. We identify two germs of Lagrangian submanifolds $L_1, L_2$ through 0
if the tangency order of $L_1$ to $L_2$ (see \cite[Def. 2.6]{DoTr}) is greater than $p$. (See also
\cite[I.1]{Je}, where the name ``contact of order'' is used.) The equivalence class is called a ``$p$-jet of a submanifold''.
In this way we obtain the {\it space of $p$-jets of
Lagrangian manifolds} denoted by ${\cal J}^p(V)$. This space is homogeneous with respect to the action of 
$\Sympl^{p}(V)$ - the {\it group of $p$-jet symplectomorphisms preserving $0\in V$}: every $p$-jet 
of a Lagrangian submanifold can be obtained from the ``distinguished'' Lagrangian submanifold $W$ 
by application of a symplectomorphism preserving $0$.

The Lagrangian Grassmannian $LG(V)$ is embedded in ${\cal J}^p(V)$ in a natural way.
On the other hand, we have the {\it Gauss map} 
$$
\pi:{\cal J}^p(V)\to LG(V)\,,
$$
which is a retraction to $LG(V)$, defined for a Lagrangian submanifold $L$ by
$
\pi(L)=T_0(L)\,, 
$
the tangent space of $L$ at $0\in L$. Denote by $\{W\}$ the point of $LG(V)$ corresponding to the linear space $W$.
The following fact and its proof stems from \cite{MPW1}.

\begin{lemma}\label{fib} The fiber of $\pi$ over $\{W\}$ is isomorphic to the linear space 
$$
\bigoplus_{i=3}^{p+1} S^i(W^*)\,.
$$
\end{lemma}
\proof The fiber $\pi^{-1}\{W\}$ consists of those (jets of)
Lagrangian submanifolds $L$ such that $T_0(L)=W$.
Every Lagrangian submanifold $L$ such that $\varrho_{|_L}$ is a submersion
is the graph of the differential of a function \ $f: W\to \C$ \
(note that $df$ acts from $W$ to $W^*$). The condition: \ $0\in L$ \
corresponds to the condition: \ $df(0)=0$, and the condition:
  \ $T_{0}(L)=W$ \ corresponds to the vanishing of the second derivatives of
$f$ at $0$. 
\qed

Thus \ $\pi: {\cal J}^p(V)\to LG(V)$ is an affine fibration. (Note that $\pi$ is not a vector bundle 
starting from $p=3$, see \cite[footnote on p. 68]{MPW1}).

Let $H$ be the subgroup of $\Sympl^p(V)$ consisting of $p$-jets of holomorphic
symplectomorphisms preserving the fibration $\varrho: V\to W$.
Two Lagrangian $p$-jets are {\it Lagrangian equivalent} if they belong to the same orbit of $H$.
{\it A Lagrange singularity class} is any closed pure dimensional algebraic subset of the manifold
$\cJ^p(V)$, which is $H$-invariant.

A Lagrange singularity class $\Sigma\subset \cJ^p(V)$ defines the class $[\Sigma]$ in the cohomology groups
\begin{equation}\label{HH}
H^*(\cJ^p(V), {\Z})\cong H^*(LG(V), {\Z})\,.
\end{equation}
This cohomology class in $H^*(LG(V), {\Z})$ will be called the {\it (Lagrangian) Thom polynomial} of $\Sigma$, and denoted $\cT^\Sigma$.

We now use Schubert calculus to investigate Lagrangian Thom polynomials, that is,
we study the expansions of Lagrangian Thom polynomials in the basis of Lagrangian Schubert classes. 
(These are the classes of the closures of the cells of a cellular decomposition of $LG(V)$, and thus they
form a basis of $H^*(LG(V), {\Z})$.) By Theorem \ref{YQ}, we have
$$
\cT^\Sigma = \sum_{\hbox{\small strict} \ {\mu}\subset \rho(n)} a_{\mu} \\Q_{\mu}(R^*)\,,
$$ 
where $a_\mu \in \Z$. 

\begin{lemma}\label{NGJ} \cite{MPW1} We have the following expression for the normal bundle of $LG(V)$ in $\cJ^p(V)$:
$$
N_{LG(V)}\cJ^p(V) \cong \bigoplus_{i=3}^{p+1} S^i(R^*)\,.
$$
\end{lemma}
In particular, $N_{LG(V)}\cJ^p(V)$ is generated by its global sections.

\begin{proposition}
Let $Z$ be a subvariety of $\cJ^p(V)$. If, using (\ref{HH}), we have
$$
[Z]=\sum a_\mu \\Q_\mu(R^*)\,,
$$
where $a_\mu \in \Z$, then all the coefficients $a_\mu$ are nonnegative.
\end{proposition}
\proof Set $G=LG(V)$, $\cJ=\cJ^p(V)$ and $N=N_G\cJ$. Denote by 
$i: G \hookrightarrow\cJ$ the inclusion. We look at the coefficients $a_{\mu}$ 
in the expression 
$$
i^*[Z]=\sum a_{\mu} \ {\\Q}_{\mu}(R^{*})= \sum a_{\mu} Y^\mu \,,
$$
where the last equality follows from Theorem \ref{YQ}. Let $Y^{\mu'}$ be the dual class to $Y^{\mu}$ (see Example \ref{LG}). 
We have \ 
$$
a_{\mu}=i^*[Z]\cdot Y^{\mu'}\,.
$$
Invoking (\ref{ha}), we may compute this last intersection number using the Chow groups of $G$.
Let $C=C_{G\cap Z} Z \subset N$ be the {\it normal cone} of $G\cap Z$ in $Z$. Denote by
$j: G\hookrightarrow N$ the zero-section inclusion. By deformation to the normal cone (see \cite[Sect. 6.1 and 6.2]{Fu1}), we have
$$
i^*[Z]=j^*[C] \ \ \ \ \ \ \ \ \hbox{(equality in $A^{*}(G)$)}\,.
$$
It follows that
$$
a_{\mu}=[C]\cdot Y^{\mu'} \ \ \ \ \ \ \ \ \hbox{(intersection in $N$)}\,.
$$

By virtue of Lemma \ref{NGJ}, the assertion now follows from Proposition \ref{CV} 
for $X=G$, $E=N$, and $Y=Y^{\mu'}$.\qed

\begin{theorem} \cite{MPW1}
For any Lagrange singularity class $\Sigma$, the Thom polynomial
$\cT^\Sigma$ is a nonnegative combination of $\\Q$-functions.
\end{theorem}

\medskip

\noindent
{\bf Question.} Does there exists a basis different (up to rescaling) from $\{Y^\mu=\\Q_\mu(R^*)\}$ with the property that any Lagrangian Thom polynomial 
has a positive expansion in that basis?

\section{Legendrian Thom polynomials}\label{tleg}

Fix $n \in \bf N$. Let $W$ be a complex vector space of dimension $n$, and let $L$ be a one dimensional complex vector space. Consider
\begin{equation}\label{V}
V:=W\oplus (W^*\otimes L)\,
\end{equation}
-- a symplectic space equipped with the {\it twisted symplectic form} $\omega\in\Lambda^2V^*\otimes L$. 

Consider a {\it contact space}
$$
V\oplus L=W\oplus (W^*\otimes L)\oplus L\,.
$$
Let $\alpha$ be a contact form on $V\oplus L$ (cf. \cite[Sect. 20.1]{AGV}).
Legendrian submanifolds of $V\oplus L$ are maximal integral
submanifolds of the form ${\alpha}$, i.e., the manifolds of dimension
$n$ with tangent spaces contained in $\Ker (\alpha)$.

To study Legendrian submanifolds (through $0$) of $V\oplus L$, we use Lagrangian submanifolds (through $0$) of $V$.
Any Legendrian submanifold in $V\oplus L$ is determined by its Lagrangian projection to $V$ 
and any Lagrangian submanifold in $V$ lifts to $V\oplus L$. 

Legendrian Thom polynomials were considered by Vassiliev \cite{Va} (see also \cite{Ka2}). In \cite{MPW2}, the space $\cJ^p(W,L)$ was constructed
(with the help of Kazarian) that can serve to address positivity questions about Legendrian Thom polynomials. This space is not a naive generalization of the space of Lagrangian $p$-jets from the previous section. Roughly speaking, one wants to parametrize the relative positions of two Lagrangian submanifolds. More precisely, we define $\cJ^p(W,L)$ to be the set of pairs of $p$-jets of Lagrangian submanifolds of $V$ 
consisting of a linear space and a submanifold whose tangent space at $0$ is $W$. For a motivation of this construction and more
details, we refer the reader to \cite[Sect. 2 and 3]{MPW2}. The projection to the first factor gives a map
$$
\pi: {\cJ}^p(W,L) \to LG(V)\,,
$$ 
which is a trivial vector bundle with the fiber equal to
$$
\bigoplus_{i=3}^{p+1} S^i(W^*)\otimes L\,.
$$
In fact, we need a relative version of this construction.
Let $X$ be a topological space, $W$ a complex rank $n$ vector bundle over $X$, and $L$ a complex
line bundle over $X$. Define a vector bundle $V$ on $X$ by (\ref{V}). Let 
$$
\tau: LG(V)\to X 
$$ 
denote the induced Lagrange Grassmann {\it bundle}. We have a {\it relative} version of the map $\pi$  
$$
\pi: {\cJ}^p(W,L) \to LG(V)\,,
$$
which is denoted by the same letter. 

The space ${\cJ}^p(W,L)$ fibers over $X$. It is equal to the pull-back
$$
{\cJ}^p(W,L) = \tau^*\left(\bigoplus_{i=3}^{p+1} S^i(W^*)\otimes L \right)\,.
$$
              
By a {\it Legendre singularity class} we mean a closed algebraic subset $\Sigma \subset {\cJ}^p(\C^n,\C)$ invariant 
with respect to holomorphic {\it contactomorphisms} of $\C^{2n+1}$. Additionally, we assume that $\Sigma$ is stable 
with respect to enlarging the dimension of $W$. Since any changes of coordinates of $W$ and $L$ induce
holomorphic contactomorphisms of $V\oplus L$, any Legendre singularity class $\Sigma$ defines
$$
\Sigma(W,L) \subset{\cJ}^p(W,L)\,.
$$   
The element $[\Sigma(W,L)]$ of $H^*({\cJ}^p(W,L),\Z)$ is called the {\it Legendrian Thom polynomial of} $\Sigma$.

In the following, we shall write ${\cJ}$ for the vector bundle  ${\cJ}^p(W,L)$.

We now use Schubert calculus to study Legendrian Thom polynomials. 
Let $L,M_1,M_2,\dots,M_n$ be one dimensional vector spaces, and let
$$
W:=\bigoplus_{i=1}^n M_i\,,\qquad V=W\oplus (W^*\otimes L)\,.
$$
We have a symplectic form $\omega$ defined on $V$ with values in $L$. The Lagrangian Grassmannian $LG(V)$ is a homogeneous 
space for the symplectic group $Sp(V) \subset \End(V)$. We fix two ``opposite'' isotropic flags $E^+$ and $E^-$ in $V$:
$$
E_j^+:=\bigoplus_{i=1}^j M_i\,,\qquad
E_j^-:=\bigoplus_{i=1}^j M_{n-i+1}^*\otimes L \,,\qquad
(j=1,2,\ldots ,n)\,.
$$
Consider two Borel groups $B^\pm\subset Sp(V)$, preserving the flags $E^\pm$.
The orbits of $B^\pm$ in $ LG(V)$ form two ``opposite'' cell decompositions 
$
\{C^\mu(E^\pm,L)\}
$ 
of $LG(V)$, labeled by strict partitions $\mu \subset \rho(n)$ (see \cite{Pr2} and Sect.\ref{Schub}). 
The cells of the $C^{-}$-decomposition are transverse to the cells of $C^{+}$-decomposition. Denote the class of the closure of $C^\mu(E^\pm,L)$ in $H^*(LG(V),\Z)$ by $Y^\mu(E^\pm,L)$. 

All these data behave functorially with respect 
to the automorphisms of the lines $L$ and $M_i$'s,
(they form a torus $(\C^*)^{n+1}$). Thus the construction of the cell decompositions can be repeated for bundles $L$ and
$\{M_i\}_{i=1}^n$ over any base $X$. We get a Lagrange Grassmann bundle 
$$
\tau: LG(V)\to X\,,
$$
endowed with two (relative) stratifications
$$
\{C^\mu(E^\pm,L) \to X\}_{\mu}\,.
$$
Suppose that $X=G/P$ is a compact manifold, homogeneous with respect to an action of a linear group $G$. Then $X$
admits a Chevalley-Bruhat cell decomposition $\{\sigma_\lambda\}$.
The subsets
$$
Z^-_{\mu \lambda}:=\tau^{-1}(\sigma_\lambda)\cap C^{\mu}(E^-,L)
$$
form an algebraic cell decomposition of $LG(V)$.
Another cell decomposition of $LG(V)$ is given by the collection of subsets
$$
Z^+_{\mu \lambda}:=\tau^{-1}(\sigma_\lambda)\cap C^{\mu}(E^+,L)\,.
$$

\begin{example}\label{Hirz} \rm
If $X=\P^1$, $W={\bf 1}$, $L=\cO(d)$ (for $d>0$),
then $LG(V)$ is the Hirzebruch surface $H_d$ which can be
presented as the sum of the space of the line bundle $L$ and
the section at infinity, \ $H_d=L \cup \P^1_\infty$.
Then $\P^1_0$, the zero section of the bundle $L$, is a stratum
of the $C^+$-decomposition and the section at infinity
$\P^1_\infty$ is a stratum of the $C^-$-decomposition. For
the cell decomposition of $X=\P^1=\C\cup \{\infty\}$, we obtain
two cell decompositions of $H_d$. Two resulting bases of cohomology 
are mutually dual with respect to the intersection product.
The closures of strata of $Z^+$-decomposition have the following property:
any effective cycle has a nonnegative intersection number with them.
This is not true for the closures of strata of $Z^-$-decomposition: for example, the
self-intersection of $\P^1_\infty$ is equal to $-d$.
\end{example}

\begin{theorem}\label{trans1} \cite{MPW2} Fix a strict partition $\mu \subset \rho(n)$ and an index $\lambda$.
Suppose that the vector bundle $\cJ$ is generated by its global sections. Then,
in $\cJ$, the intersection of $\Sigma(W,L)$ with the closure of
any $\pi^{-1}({Z^-_{\mu \lambda}})$ is represented by a nonnegative cycle.
\end{theorem}
The proof in \cite{MPW2} is based (apart from the Schubert calculus for $LG(V)\to X$) on some variant of the Bertini-Kleiman theorem.

\medskip

We apply the theorem in the situation when all $M_i$ are equal to the same line bundle $M$ 
(and then $W=M^{\oplus n}$) and $M^{-m} \otimes L$ is generated by its global sections for $m\geq 3$.
             
Consider the following three cases: the base is always $X={\bf P}^n$ and

\begin{itemize}
\item $L_1=\cO(-2)$,  \ $M_1=\cO(-1)$, or 

\item $L_2=\cO(1)$,  \  $M_2={\bf 1}$, or

\item $L_3=\cO(-3)$, \ $M_3=\cO(-1)$.
\end{itemize}

We obtain the symplectic bundles $V_i=M_i^{\oplus n}\oplus(M_i^*\otimes L_i)^{\oplus n}$ with twisted symplectic forms $\omega_i$
for $i=1,2,3$.

These three cases were crucial to discover and prove the forthcoming Theorem \ref{main2}.
Case 1 was the subject of \cite[Rem. 14]{MPW1}. 
In Case 2, the integral cohomology $H^*(LG(V),\Z)$ is
isomorphic to the ring of Legendrian characteristic classes up to degree $n$;
the $Z^-$-decomposition of $LG(V)$ gives us another basis of cohomology.
In Case 3, the cohomology of $LG(V)$ is isomorphic, up to degree $n$, to the ring of Legendrian 
characteristic classes, provided we invert the number 3 this time.
The positivity property in Case 1 was known ({\it loc.cit.}), whereas in Cases 2 and 3, it was Kazarian
who suggested the positivity.
       
To overlap all these three cases we consider 
$X:={\bf P}^n\times {\bf P}^n$
and set
$$
W:= p_1^* \cO(-1)^{\oplus n}\,,\qquad L:=p_1^*\cO(-3)\otimes p_2^*\cO(1)\,,
$$
where $p_i:X\to {\bf P}^n$, $i=1,2$, are the projections. Restricting the bundles $W$ and $L$ to the diagonal
or to the factors, we obtain the three cases considered above.

The space $LG(V)$ has a cell decomposition
$Z^+_{{\mu} \lambda}$, where $\mu$ runs over strict partitions contained
in $\rho(n)$, and $\lambda =(a,b)$ with $a$ and $b$ natural numbers
smaller than or equal to $n$. The classes of closures of the cells
of this decomposition give a basis of the cohomology of $LG(V)$.

Let $v_1$ and $v_2$ be the first Chern classes of $p_1^*(\cO(1))$ and $p_2^*(\cO(1))$.
We have
\begin{equation}\label{factor}
[\overline{Z^+_{\mu,a,b}}]=Y^\mu(E^+,L) \ v_1^a v_2^b\,.
\end{equation}

\begin{theorem}\label{main2} \cite{MPW2} Let $\Sigma$ be a Legendre singularity class.
Then $[\Sigma(W,L)]$ has nonnegative coefficients in the basis $\{[\overline{Z^+_{\mu,a,b}}]\}$.
\end{theorem}

\begin{example} \rm Using the names of singularities from \cite{Ka2}, we display some Legendrian Thom polynomials in the basis 
from the theorem. The bold terms give the Thom polynomials of the corresponding Lagrange singularities.

\medskip

\noindent 
$A_2$: $\bf\\Q_1$

\noindent 
$A_3$: ${\bf 3\\Q_2}+v_2\\Q_1$

\noindent 
$A_4$: ${\bf 12\\Q_3+3\\Q_{21}}+(3v_1+7v_2)\\Q_2+(v_1v_2+v_2^2)\\Q_1$

\noindent 
$D_4$: $\bf\\Q_{21}$

\noindent
$P_8$: $\bf\\Q_{321}$.

\noindent 
$A_5$: ${\bf     60\\Q_4+27\\Q_{31}}+
      (6v_1+16v_2)\\Q_{21}+
           (39v_1+47v_2)\\Q_3+$

               $ (6v_1^2+22v_1v_2+12v_2^2)\\Q_2+
(2v_1^2v_2+3v_1v_2^2+v_2^3)\\Q_1$

\noindent 
$D_5$: ${\bf 6\\Q_{31}}+4v_2\\Q_{21}$,

\noindent 
$P_9$: ${\bf 12\\Q_{421}}+12v_2\\Q_{321}$.

\end{example}

Using the theorem, one constructs, in the ring of Legendrian characteristic classes, a one-parameter family of bases such 
that any Legendrian Thom polynomial has, in any basis from the family, an expansion with nonnegative coefficients (see \cite{MPW2}). 

\begin{remark} \rm
Positive descriptions of the coefficients of Schur function expansions of Thom polynomials are known for several series
of singularity classes of maps, see a survey article \cite{OzPr}. For those coefficients, which do not admit such descriptions,
it is interesting to establish their bounds. For the Legendrian Thom polynomials, this issue is discussed in \cite[Sect. 9 and 10]{MPW2}.
In particular, one examines there how positivity of Thom polynomials of maps to curves implies some upper bounds on the coefficients
of Legendrian Thom polynomials.
\end{remark}

\begin{remark} \rm
It is shown in \cite[Sect. 10]{MPW2} that the Thom polynomials of nonempty stable Lagrangian and Legendrian singularity classes
are nonzero.
\end{remark}

%\begin{note} \rm Huh \cite{Hu} has recently shown the validity of the positivity conjecture, stated in \cite{AlMi}, about 
%the Schwartz-MacPherson classes of the Schubert cells in Grassmannians.
%\end{note} 

\end{document}